\documentclass[a4paper]{elsart3-1}
\usepackage{amsmath,amssymb}
\usepackage{pstricks,pst-node}
\usepackage{enumerate}

\usepackage[francais,english]{babel}
\usepackage[T1]{fontenc}
\usepackage[latin1]{inputenc}

\newtheorem{theorem}{Theorem}[section]

\newtheorem{e-proposition}[theorem]{Proposition}

\newtheorem{e-definition}[theorem]{Definition\rm}

\newtheorem{theoreme}{Th\'eor\`eme}[section]
\newtheorem{lemme}[theoreme]{Lemme}
\newtheorem{proposition}[theoreme]{Proposition}

\renewcommand{\theequation}{\arabic{equation}}
\def\og{\leavevmode\raise.3ex\hbox{$\scriptscriptstyle\langle\!\langle$~}}
\def\fg{\leavevmode\raise.3ex\hbox{~$\!\scriptscriptstyle\,\rangle\!\rangle$}}

\journal{the Acad\'emie des sciences}
\begin{document}
\centerline{Combinatoire}
\begin{frontmatter}

\selectlanguage{francais}
\title{Inversion dans les tournois \thanksref{label1}}
\thanks[label1]{Les auteurs d\'edient ce texte \`a la m\'emoire de Roland Fra\"{i}ss\'e.}

\author[authorlabel1]{Houmem BELKHECHINE},
\ead{houmem@gmail.com}
\author[authorlabel2]{Moncef BOUAZIZ},
\ead{moncef.bouaziz@laposte.net}
\author[authorlabel3]{Imed BOUDABBOUS},
\ead{imed.boudabbous@fsegs.rnu.tn}
\author[authorlabel4]{Maurice POUZET}
\ead{maurice.pouzet@univ-lyon1.fr}
\address[authorlabel1]{Facult\'e des Sciences de Gab\`es, Gab\`es, Tunisie}
\address[authorlabel2]{Institut des technologies m\'edicales, Tunis, Tunisie}
\address[authorlabel3]{Institut Pr\'eparatoire aux \'Etudes d'Ing\'enieurs de Sfax, Sfax, Tunisie}
\address[authorlabel4]{ICJ, Math\'ematiques, Universit\'e Claude Bernard Lyon 1, Lyon, France et Department of Mathematics and Statistics, The University of Calgary, Calgary, Alberta, Canada}


\begin{abstract}
\selectlanguage{francais}
Nous consid\'erons  la transformation qui inverse tous les arcs d'une partie  $X$ de l'ensemble des sommets d'un tournoi
$T$.  L'\emph{indice} de $T$, not\'e $i(T)$, est
le plus petit nombre de parties dont il faut inverser les arcs pour ramener $T$
\`a un tournoi acyclique. Il appara\^it que les tournois critiques et les tournois $(-1)$-critiques peuvent \^etre d\'efinis au moyen d'inversions, les premiers \'etant d'indice un ou deux, les seconds d'indice au plus quatre. On peut voir  $i(T)$ comme le minimum de la distance de $T$ aux tournois acycliques d\'efinis sur le m\^eme ensemble de sommets;  la distance entre deux tournois $T$ et $T'$ peut \^etre \'egalement interpr\'et\'ee comme la  \emph{dimension bool\'eenne} d'un graphe, celui-ci \'etant la  somme bool\'eenne de $T$ et $T'$. Sur  $n$ sommets, la distance maximale vaut $n-1$ tandis que   $i(n)$,  le maximum des indices des tournois \`a $n$ sommets,
 satisfait les in\'egalit\'es  $\frac{n-1}{2} - \log_{2}n\leq i(n) \leq n-3$ pour  $n \geq 4$. Soit $\mathcal {I}_{m}^{<\omega}$
(resp. $\mathcal {I}_{m}^{\leq\omega}$), la classe des tournois finis (resp.
au plus d\'enombrables) $T$ tels que $i(T) \leq m$. La classe $\mathcal {I}_{m}^{<\omega}$ est d\'etermin\'ee par un nombre fini d'obstructions; nous  donnons une description
morphologique des \'el\'ements de $\mathcal{I}_{1}^{<\omega}$ et d\'ecrivons ses obstructions. Nous d\'ecrivons aussi un tournoi
universel de la classe $\mathcal {I}_{m}^{\leq\omega}$.
{\it Pour citer cet article~: Houmem Belkhechine, Moncef Bouaziz, Imed Boudabbous, Maurice Pouzet, C. R.
Acad. Sci. Paris (2010).}
\vskip 0.5\baselineskip

\selectlanguage{english} \noindent{\bf Abstract} \vskip
0.5\baselineskip \noindent {\bf Inversions in tournaments.}
We consider  the transformation reversing  all arcs  of a subset $X$ of the vertex set of a tournament $T$. The \emph{index} of $T$, denoted by
$i(T)$, is the smallest number of subsets that must be reversed to make $T$ acyclic. It turns out that
critical tournaments and $(-1)$-critical tournaments  can be defined in terms of inversions (at most two for the former, at most four for the latter). We interpret   $i(T)$ as the minimum distance of  $T$ to the transitive tournaments on the same vertex set, and we interpret the distance between two tournaments $T$ and  $T'$ as the \emph{Boolean dimension} of a graph, namely the Boolean sum of  $T$ and $T'$. On $n$ vertices, the maximum distance is at most $n-1$, whereas  $i(n)$,  the maximum of $i(T)$ over the tournaments
on $n$ vertices, satisfies  $\frac {n-1}{2} - \log_{2}n
\leq i(n) \leq n-3$, for  $n \geq 4$. Let $ \mathcal{I}_{m}^{< \omega}$ (resp. $\mathcal{I}_{m}^{\leq \omega}$) be the class of finite (resp.
at most countable) tournaments $T$ such that $i(T) \leq m$. The class $\mathcal {I}_{m}^{< \omega}$ is determined by finitely many obstructions. We give a morphological description
of the members of $\mathcal {I}_{1}^{< \omega}$ and a description of the critical obstructions. We give an explicit description of  an universal tournament of the class $\mathcal{I}_{m}^{\leq \omega}$.
{\it To cite this article: Houmem Belkhechine, Moncef Bouaziz, Imed Boudabbous, Maurice Pouzet, C. R.
Acad. Sci. Paris (2010).}
\end{abstract}
\end{frontmatter}

\selectlanguage{english}
\section*{Abridged English version}
Let $T$ be a tournament.  Let  $V(T)$ be its vertex set and $A(T)$ be its arc set. An {\it inversion} of an arc $a := (x,y)\in A(T)$   consists  to replace the arc $a$ by
$a^{\star} := (y,x)$ in $A(T)$. For a subset $X
\subseteq V(T)$, let  $Inv(T, X)$ be the tournament  obtained from
$T$ after reversing all arcs $(x,y) \in
A(T)\cap (X\times X)$. For example, $Inv(T, V(T))$ is $T^{\star}$, the {\it dual}  of $T$. For a finite sequence $(X_{i})_{ i
< m}$ of subsets of $V(T)$,  let $Inv( T, (X_{i})_{i< m})$ be the tournament obtained from $T$
by reversing successively all the arcs in each of the  subsets $X_i$, $i<m$,  that is the
tournament equal to $T$ if $m=0$ and to $Inv(Inv(T, (X_{i})_{i< m-1}), X_{m-1})$ if $m\geq 1$. The {\it inversion index} of $T$, denoted by $i(T)$,  is  the least integer
$m$  such that there is a sequence $(X_{i})_{i< m}$ of subsets of $V(T)$ for which $Inv( T, (X_{i})_{i< m})$ is
acyclic. This is a variant of the {\it Slater index} of $T$ (the minimum number
of arcs which should be reversed to make it
acyclic, \cite{SL}). Our motivation originates in the study of critical tournaments. Indeed, the critical tournaments characterized in \cite {ST} can be easily defined from acyclic tournaments by means of one or two inversions  whereas the $(-1)$-critical tournaments, characterized  in \cite {HIJ2}, can be defined by means of two, three or  four inversions \cite {belkhechine-these};   an other interest comes from the point of view of logic. We present some general properties of the inversion index and of the class  $\mathcal I_{m}$ of tournaments $T$ having inversion index at most $m$, with a particular emphasis on the subclasses $ \mathcal {I}_{m}^{< \omega}$ and $ \mathcal {I}_{m}^{\leq \omega}$ made respectively of finite and at most countable members of $\mathcal I_{m}$.  Part of these results are included in \cite{belkhechine-these}.  We use  tools from the theory of relations in the vein of Fra\"{\i}ss\'e, referring to \cite{F}  for the notions of relational structure, embeddability, classes closed under embeddability, that we call here \emph{hereditary classes}, bounds, age and free operator. We leave open the road for algorithmic considerations.

Let   $\mathcal{T}_{V}$ be  the set of tournaments $T$ on a fixed set $V$ of vertices. Pairs  $(T,T')$ of distinct members of  $\mathcal{T}_{V}$ such that $T'=Inv(T, X)$ for some $X\subseteq V$ form the edges of a irreflexive and symmetric  graph on $\mathcal{T}_{V}$. With respect to  the graphic distance associated with this graph, the inversion index of $T\in \mathcal{T}_{V}$ is then the  minimum distance of $T$ to the  acyclic members of $\mathcal{T}_{V}$. We define the  \emph{Boolean dimension} of  a graph $G$ (irreflexive and symmetric) as the least integer $m$ such that $G$ can be represented by the non orthogonality relation on  the vector space $\mathbb F_2^{m}$ equipped with the ordinary scalar product.
\begin{theorem} \label{THM distanc}
The graphic distance $d(T,T')$ between  two members $T,T'$ of $\mathcal{T}_{V}$ is the Boolean dimension of the Boolean sum $T\dot{+}T'$. If $V$ has $n$ elements, this distance is at most $n-1$. It is attained if $T'= T\dot{+}P$ where $P$ is any path on $V$.

\end{theorem}
For $n \in \mathbb{N}$, let $i(n)$ be  the maximum of the inversion index of tournaments on $n$ vertices.

\begin{theorem} $\frac{n-1}{2} - \log_{2} n
\leq i(n) \leq n-3$ for all integer $n \geq 4$.
\end{theorem}

If $T\in \mathcal T_{V}$ and $(X_{i})_{i< m}$ is  a sequence of subsets of $V$, we consider  the pair $(T, (X_{i})_{i< m})$ as a relational structure made of the set $V$, the binary relation $A(T)$ and the unary relations $X_i$ for $i<m$. Let $\mathcal O_{m}$, resp. $\mathcal O_{m}^{\leq \omega}$, resp. $\mathcal O_{m} ^{< \omega}$,  be the class of these $(T, (X_{i})_{i< m})$ where $T$ is acyclic and its size is arbitrary, resp. at most countable, resp. finite. The transformation of  each $(T, (X_{i})_{i< m})$ into $Inv(T, (X_{i})_{i< m})$  defines a free operator from   $\mathcal O_{m}$ onto $ \mathcal {I}_{m}$. With this formalism follows readily that \emph{a tournament $T$ belongs to  $ \mathcal {I}_{m}$  if and only if every finite subtournament of $T$ belongs to  $\mathcal {I}_{m}$}.  The class $\mathcal O^{<\omega}_{m}$ is  a  Fra\"{\i}ss\'e class (it is hereditary and has the amalgamation property) hence  $\mathcal O^{\leq \omega}_{m}$ contains an homogenous structure with age  $\mathcal O^{<\omega}_{m}$. This structure is unique up to isomorphisms. We denote it by $C(m)$. Let $W(m):=Inv(C(m))$. This tournament is universal for $\mathcal {I}_{m}^{\leq \omega}$, that is belongs  to $\mathcal{I}_{m}^{\leq \omega}$ and embeds all members of  $\mathcal{I}_{m}^{\leq \omega}$. For $m = 0$, $W(0):= \underline{\mathbb{Q}}$. We give an explicit description of $W(m)$ for all others values of $m \geq 1$.
If $B$ is a class of finite tournaments, we denote by $Forb_{\mathcal T}(B)$ the class of finite tournaments in which no member of  $B$ is embeddable.
If $\mathcal{C}$ is a hereditary class  of finite tournaments, a \emph{bound} of  $\mathcal{C}$ in the class $\mathcal T$ of tournaments is a finite tournament $T$ not belonging to $\mathcal{C}$ such that for all $x \in V(T)$, the tournament $T-x$ belongs to $\mathcal{C}$. We
denote by $B_{\mathcal T}(\mathcal{C})$ the collection of these bounds. As it is well known for arbitrary hereditary classes of finite relational structures, a hereditary class of finite tournaments is determined by its bounds; in fact  $\mathcal{C} = Forb_{\mathcal T}(B_{\mathcal T}(\mathcal {C}))$. The test given in \cite {pouzet2} and Higman theorem on words yield:

\begin{theorem} The  class $\mathcal {I}_{m}^{< \omega}$
has only finitely many bounds, these bounds being considered up to isomorphisms.
\end{theorem}
For example,  the $3$-cycle $C_{3}$ is, up to isomorphisms,  the unique bound of the class $\mathcal {I}_{0}^{< \omega}$. Let $n\in \mathbb N$. We set $\mathbb {N}_{<n}:=\{i\in \mathbb {N}: i<n\}$. We denote by $\underline{n}$ the tournament whose vertex set is $\mathbb {N}_{<n}$  and whose arcs are pairs $(i, j)$ such that $0\leq i<j<n$.

\begin{theorem} \label{THM fini de bornes2}
The bounds of the class $\mathcal{I}_{1}^{< \omega}$ are, up to isomorphisms,
$B_{6} := Inv (\underline{6}, (\{0,3,4\},\{1,4,5\}))$, $C_{3}. \underline{2} := Inv (\underline{6}, (\{0,2\},\{3,5\}))$,
$D_{5} := Inv (\underline{5}, (\{1,3\},\{0,4\}))$, $T_{5} := Inv(\underline{5}, (\{0,2,4\}, \{1,3\}))$
and $V_{5} := Inv(\underline{5}, (\{0,4\}, \{2,4\}))$. In particular, $\mathcal{I}_{1}^{<\omega} = Forb_{\mathcal T}(\{B_{6}, C_{3}.\underline {2}, D_{5}, T_{5}, V_{5}\})$.
\end{theorem}

Let $(T_{i})_{i \in V}$  be a family of tournaments whose vertex sets $V(T_{i})$ are pairwise disjoint. If  $T$ is  a tournament with vertex set $V$,  we denote by $\sum_{i \in T} T_{i}$
the  \emph{lexicographical sum of the $T_{i}$ 's indexed by $T$}. When $V = \mathbb{N}_{<n}$, where $n \in \mathbb{N} \backslash \{0\}$,
$\sum_{i \in T} T_{i}$ is also denoted by $T(T_{0}, \cdots, T_{n-1})$. It turns  out that $i(\sum_{j \in T} T_{j})=i(T)$ provided that the $T_j$' s
are acyclic and non empty. A tournament $T$ is \emph{acyclically indecomposable} if no acyclic autonomous subset of $T$ has more than one element  \cite{C-J}. Since every tournament is a lexicographical sum of acyclic tournaments indexed by some  acyclically indecomposable tournament $T$ \cite{B-P},  \emph{the members of $\mathcal{I}_{m}$ are the lexicographical sums of  acyclic tournaments indexed  by acyclically indecomposable members  of $\mathcal {I}_{m}$}. Alternatively, {\it  the bounds of  $\mathcal{I}_{m}^{< \omega}$ in $\mathcal T$ are acyclically indecomposable}.
\begin{theorem}
A tournament $T$ with $\mid V(T)\mid \geq 2$ is an acyclically indecomposable member of $\mathcal{I}_{1}^{<\omega}$
if and only if $T$ is isomorphic to $U_{2n+1}$, $\underline{2}(\underline{1}, U_{2n+1})$, $\underline{2}(U_{2n+1}, \underline{1})$ or $\underline{3}(\underline{1}, U_{2n+1}, \underline{1})$, where $n \geq 1$ and $U_{2n+1}: = Inv(\underline{2n+1}, 2\mathbb{N}_{<n+1})$.
\end{theorem}

\selectlanguage{francais}

\selectlanguage{francais}

\section{Terminologie}
Dans cette note, nous consid\'erons essentiellement des tournois et des graphes. Nous consid\'erons ceux-ci comme des digraphes sans boucles. Nous rappelons quelques notions concernant les digraphes et renvoyons \`a \cite {Bo} et \cite{F} pour la terminologie non d\'efinie. Nous rappelons qu'un \emph{digraphe} ou \emph {graphe dirig\'e} est un couple $D$ form\'e d'un ensemble $V$, dont les \'el\'ements sont les \emph{sommets} de $D$, qu'on note $V(D)$,  et d'une partie,  not\'ee $A(D)$,  du produit $V\times V$ dont les \'el\'ements sont les \emph{arcs} de $D$. Le \emph{dual} de $D$, not\'e $D^{\star}$, est le digraphe ayant m\^emes sommets et pour arcs les couples $(x,y)$ tels que $(y,x)\in A(D)$. Si $X$ est une partie de $V(D)$ alors $D_{\restriction X}:=(X,  A(D) \cap (X \times X)) $ est le \emph{digraphe induit par $D$ sur  $X$}. Si $x\in V(D)$ nous notons  $D-x$ le digraphe induit sur $V(D) \setminus \{x\}$. Un digraphe $D$ {\it s'abrite} (ou se plonge) dans un
digraphe $D'$ lorsque $D$ est isomorphe au digraphe induit par  $D'$ sur une de ses parties. Ainsi un tournoi \emph{acyclique} ou  {\it transitif}
est un tournoi qui n'abrite pas le  $3$-cycle $C_{3} := (\{0,1,2\},\{(0,1), (1,2), (2,0)\})$. Nous d\'esignons par $\mathcal D_V$ l'ensemble des digraphes ayant $V$ comme ensemble de sommets. Soit $D\in \mathcal D_V$. Pour $x,y\in V$ nous posons $D(x,y)=1$ si $(x,y)\in A(T)$ et sinon $D(x,y)=0$. Une partie $X$ de $V$ est
un \emph{intervalle} de $D$  lorsque $D(x, y)=D(x', y)$ et $D(y, x)=D(y,x')$ pour tout $y\in V\setminus X$ et $x,x'\in X$.  Par exemple,
 $\emptyset $, $V$, et $\{x\}$ o\`u $x \in V$ sont des intervalles de $D$ dits \emph{triviaux}. Le digraphe $D$   est {\it ind\'ecomposable}  si tous ses intervalles sont triviaux et  {\it d\'ecomposable} dans le cas contraire.
Si  $(D_{i})_{i \in V}$  est une famille de digraphes dont les ensembles de sommets  $V_{i}:= V(D_{i})$ sont deux \`a deux   disjoints et si  $D$ est un digraphe  ayant  $V$ comme ensemble de sommets, la \emph{somme lexicographique des $D_i$ index\'ee par $D$} est le digraphe not\'e $\sum_{i \in D} D_{i}$ dont les sommets sont les \'el\'ements de  $\bigcup\limits_{i\in V}V_{i}$ et  les arcs les couples $(x,y)$  tels que ou bien $(x,y)$ est
un arc de l'un des $D_{i}$, ou bien $x \in V_{i}$, $y \in V_{j}$,
avec $i \neq j \in V$, et $(i,j)$ est un arc de $D$. Pour $n\in \mathbb N$, notons $\mathbb{N}_{<n}:=\{i\in \mathbb {N}: i<n\}$,
$2 \mathbb{N}_{<n} := \{2i : i \in \mathbb{N}_{<n}\}$, $2 \mathbb{N}_{<n} + 1 := \{2i+1 : i \in \mathbb{N}_{<n}\}$, et d\'esignons par  $\underline{n}$ le tournoi ayant $\mathbb {N}_{<n}$ comme ensemble de sommets et pour arcs  les  couples  $(i, j)$ tels que $0\leq i<j<n$. Lorsque $V = \mathbb{N}_{<n}$, $n \in \mathbb{N} \backslash \{0\}$, la somme  $\sum_{i \in D} D_{i}$ est aussi not\'ee $D(D_{0}, \cdots, D_{n-1})$ et $\underline n(D_{0}, \cdots, D_{n-1})$ lorsque $D=\underline n$.

\section{Inversion et indice d'inversion}

Soit $T$ un tournoi. Une inversion d'un arc $a := (x,y)\in A(T)$ dans $T$ consiste \`a remplacer  l'arc $a$ par
$a^{\star} := (y,x)$. Pour $X \subseteq V(T)$, nous notons $Inv(T, X)$ le tournoi obtenu en inversant tous les arcs $a\in A(T)\cap (X\times X)$. Par exemple, lorsque $X = V(T)$, $Inv(T,X) =T^{\star}$. Si $(X_{i})_{ i
< m}$ est une suite finie de parties de $V(T)$, nous notons $Inv( T, (X_{i})_{i< m})$ le tournoi obtenu \`a partir de  $T$
en inversant successivement tous les  arcs ayant leurs sommets dans $X_i$, ceci pour $i<m$. Autrement dit $Inv(
T, (X_{i})_{i< m})=T$ si $m=0$ et $Inv(Inv(T, (X_{i})_{i< m-1}), X_{m-1})$ si $m\geq 1$. De fa\c{c}on \'equivalente, un arc $(x,y)\in A(T)$ est invers\'e si et seulement si le nombre d'indices $i$ tels que $\{x,y\} \subseteq X_i$ est impair.
Cette notion d'inversion donne lieu \`a une pr\'esentation simple des tournois critiques et des tournois $(-1)$-critiques. Rappelons qu'un sommet $x$ d'un tournoi non vide, fini et ind\'ecomposable $T$ est dit {\it critique} si le tournoi
$T-x$ est d\'ecomposable et que $T$ est dit \emph{critique}, resp. \emph{$(-1)$-critique},  si tous ses sommets sont critiques, resp. si un et un seul sommet de $T$ est non critique. Les tournois critiques  ont \'et\'e introduits et caract\'eris\'es en \cite{ST}, les (-1)-critiques en
\cite{HIJ2}.  Avec la notion d'inversion, ils peuvent \^etre  d\'ecrits comme suit \cite{belkhechine-these}.

\begin{theoreme}
\begin{sloppypar}
\`A des isomorphismes pr\`es,  les tournois critiques sont les tournois $U_{2n+1} := Inv(\underline{2n+1}, 2 \mathbb{N}_{<n+1})$,
 $T_{2n+1} := Inv(\underline{2n+1}, (2 \mathbb{N}_{<n+1}, 2 \mathbb{N}_{<n} +1) )$ et $V_{2n+1} := Inv (\underline{2n+1}, (2 \mathbb{N}_{<n+1}, 2 \mathbb{N}_{<n}))$, o\`u $n \geq 2$.
 \end{sloppypar}
 \end{theoreme}
 \begin{theoreme}
 \begin{sloppypar}

\`A des isomorphismes pr\`es, les tournois (-1)-critiques sont les tournois
$E_{2n+1}^{2k+1} := Inv(\underline{2n+1}, (2\mathbb{N}_{<n+1}, 2\mathbb{N}_{<k+1}, 2\mathbb{N}_{<n+1}\setminus 2\mathbb{N}_{<k+1}))$, $F_{2n+1}^{2k+1} := Inv(\underline {2n+1},  (2\mathbb{N}_{<n+1}, 2\mathbb{N}_{<n+1}\setminus 2\mathbb{N}_{<k+1}))$, $G_{2n+1}^{2k+1} := Inv(\underline {2n+1}, (2\mathbb{N}_{<n+1}, 2\mathbb{N}_{<n}, 2\mathbb{N}_{<k+1}))$, $H_{2n+1}^{2k+1} := Inv(\underline {2n+1}, (2\mathbb{N}_{<k+1}, 2\mathbb{N}_{<k}, 2\mathbb{N}_{<n+1} \setminus 2\mathbb{N}_{<k}, 2\mathbb{N}_{<n} \setminus 2\mathbb{N}_{<k}))$, $(F_{2n+1}^{2k+1})^{\star}$ et $(G_{2n+1}^{2k+1})^{\star}$, o\`u $n \geq 3$ et $k \in \{1, \cdots, n-2\}$.
\end{sloppypar}

\end{theoreme}

Nous d\'efinissons l'{\it indice d'inversion} d'un
tournoi  $T$, not\'e $i(T)$,  comme le plus petit entier
$m$, s'il existe, pour lequel  $Inv(T, (X_{i})_{i <m})$ est un tournoi acyclique, sinon $i(T)$ est infini.  Nous pr\'esentons quelques r\'esultats simples concernant cette notion.

\section{Dimension bool\'eenne des graphes, distance et indice d'inversion}
 Soient  $D, D'\in \mathcal D_V$;  la \emph{somme bool\'eenne de $D$ et $D'$} est le digraphe not\'e $D\dot{+} D'\in \mathcal D_V$ dont l'ensemble des arcs est la diff\'erence sym\'etrique  $A(D)\Delta A(D')$; autrement dit $(D\dot{+}D')(x,y)=D(x,y) + D'(x,y)$ o\`u la somme est prise modulo $2$. Ceci permet de voir  $\mathcal D_V$ comme espace vectoriel sur le corps $\mathbb F_2$. Le sous-ensemble $\mathcal G_V$ form\'e des graphes (sans boucles et sym\'etriques) est un sous-espace de $\mathcal D_V$ et le sous-ensemble  $\mathcal T_V$ des tournois est un translat\'e de $\mathcal G_V$. Soit $G\in \mathcal G_V$. Si  $F$ est un espace vectoriel sur le corps $\mathbb F_2$ et $\varphi$ une forme bilin\'eaire et sym\'etrique sur $F$, une  \emph{repr\'esentation} de $G$ dans $(F,\varphi)$ est une application $f$ de $V$ dans $F$ telle que $G(x,y)=\varphi (f(x),f(y))$ pour tous $x, y\in V$ tels que $x\not =y$. On observera que si $f$ est une repr\'esentation alors pour tout  $v\in F$, $f^{-1}(v)$ est un intervalle de $G$;  s'il a au moins deux \'el\'ements, c'est un stable ou une clique de $G$ suivant que $v$ est isotrope ou non.  Noter que $G$  a une repr\'esentation dans l'espace de dimension $0$, resp. $1$, si et seulement si $G$ est un stable, resp. est la somme directe d'une clique et d'un stable.  La notion de repr\'esentation donne lieu \`a trois notions de dimension (binaire, isotropique ou bool\'eenne) suivant la nature de $\varphi$. La \emph{dimension bool\'eenne} de $G$ est le plus petit entier $m$ tel que $G$ admette une repr\'esentation dans  $(F, \varphi)$, o\`u $F= (\mathbb F_2)^m$ et $\varphi(u,v)=\sum_{i<m}u_iv_i$ modulo $2$.
\begin{proposition} La dimension bool\'eenne d'un graphe $G$ est le  plus petit nombre $m$ ($m\in \mathbb N$) de parties $X_i$, $i<m$,  de $V(G)$ pour lesquelles une paire $e:=\{x,y\}$ d'\'el\'ements distincts de $V(G)$ est une ar\^ete de $G$ si et seulement si le nombre de  parties $X_i$  contenant la paire $e$ est impair. Si $G$ a $n$ sommets,  sa dimension bool\'eenne est au plus $n-1$. Cette valeur maximum est atteinte par n'importe quel chemin. \end{proposition}


Les couples  $(T,T')$ d'\'el\'ements distincts de $\mathcal{T}_{V}$ tels que  $T'=Inv(T, X)$ pour une partie $X\subseteq V$ forment les arcs d'un graphe sans boucle et sym\'etrique sur  $\mathcal{T}_{V}$.

\begin{theoreme} La distance graphique $d(T,T')$  entre deux tournois $T,T'\in \mathcal{T}_{V}$ est la dimension bool\'eenne  de la somme bool\'eenne $T\dot{+}T'$. Si $V$ a $n$ \'el\'ements, cette distance est au plus $n-1$. Elle est atteinte si $T'= T\dot{+}P$ o\`u $P$ est n'importe quel chemin. L'indice d'inversion de $T\in \mathcal{T}_{V}$ est le minimum de la distance de  $T$ aux tournois acycliques appartenant \`a $\mathcal{T}_{V}$.
\end{theoreme}

\begin{lemme}\label{sum}  $i(\sum_{j \in T} T_{j})=i(T)$ pour toute famille de tournois acycliques non vides index\'ee par un tournoi.
\end{lemme}

Pour $n \in \mathbb{N}$, nous  d\'esignons par $i(n)$  l'indice d'inversion maximum des tournois \`a $n$ sommets.
Il satisfait l'in\'egalit\'e $i(n)\leq i(n-1)+1$. Pour $n \geq 4$, on obtient $i(n) \leq n-3$. Pour  $N \in \mathbb{N}$,
$\mid \{T \in \mathcal{T}_{\mathbb{N}_{<n}} : i(T) < N \} \mid \leq n! 2^{n(N-1)}$, donc pour  $m \in \mathbb{N}$ tel que $2^{\frac{m(m-1)}{2}} > m! 2^{m(N-1)}$, il existe
un tournoi $T$ d'ordre $m$ tel que $i(T) \geq N$. Ainsi:
\begin{theoreme} \label{THM 0}
Pour tout entier $n \geq 4$, $\frac{n-1}{2} - \log_{2} n
\leq i(n) \leq n-3$.
\end{theoreme}

\section{Classes de tournois d'indice born\'e}
Soit $\mathcal T$ la classe des tournois. Pour $m \in \mathbb{N}$,  nous notons $\mathcal{ I}_m$ la classe des tournois d'indice au plus $m$ et $\mathcal {I}_{m}^{<\omega}$, resp. $\mathcal {I}_{m}^{\leq\omega}$, la sous-classe de ceux qui sont finis, resp. au
plus d\'enombrables. Nous \'etudions ces classes au moyen de concepts  de la th\'eorie des relations. Nous nous r\'ef\'erons  \`a \cite{F} pour les notions concernant les structures relationnelles, e.g. abritement, classes closes pour l'abritement, que nous appelons ici
 \emph{classes h\'er\'editaires}, bornes, \^age, interpr\'etabilit\'e libre et op\'erateur libre. Si $T$ est un tournoi et  $(X_{i})_{i< m}$ une suite de  parties de $V(T)$, nous consid\'erons le couple  $(T, (X_{i})_{i< m})$ comme une structure relationnelle faite de l'ensemble $V(T)$, de la relation binaire $A(T)$ et des relations unaires $X_i$, $i<m$. Soit $\mathcal O_{m}$, resp. $\mathcal O_{m}^{\leq \omega}$, resp. $\mathcal O_{m} ^{< \omega}$, la classe des structures $(T, (X_{i})_{i< m})$ dans lesquelles  $T$  est acyclique et le cardinal de $V(T)$ est arbitraire, resp. au plus d\'enombrable,  resp. finie.  Soit $\mathcal F_{m}$ la classe des tournois librement interpr\'etables par un \'el\'ement de $\mathcal O_{m}$. Cette classe contient $\mathcal I_{m}$; en fait la transformation de chaque  $(T, (X_{i})_{i< m})\in\mathcal O_{m}$ en  $Inv(T, (X_{i})_{i< m})$  est un op\'erateur libre qui transforme $\mathcal O_{m}$ en $\mathcal I_{m}$. Par compacit\'e,   \emph{un tournoi  $T$ appartient \`a  $ \mathcal {I}_{m}$ si et seulement si  il en va de m\^eme de tout sous-tournoi fini de $T$}. %

 Une classe $\mathcal{C}$ de tournois est \emph{h\'er\'editaire} si tout tournoi
qui s'abrite dans un tournoi de $\mathcal{C}$ est encore dans  $\mathcal{C}$. Si $B$ est une classe de tournois finis,
nous  d\'esignons  par $Forb_{\mathcal T}(B)$ la classe des
tournois finis n'abritant aucun \'el\'ement de  $B$. C'est une classe h\'er\'editaire;  en fait, toute classe
h\'er\'editaire de tournois finis peut \^etre obtenue de cette fa\c{c}on. Si  $\mathcal{C}$ est une telle classe, une \emph{borne} de $\mathcal C$ dans $\mathcal T$ est tout tournoi fini $T$ minimal pour l'abritement \`a ne pas \^etre dans $\mathcal C$. Autrement dit,
 $T\not \in \mathcal{C}$ et  $T-x\in \mathcal{C}$ pour tout $x \in V(T)$. Ainsi, si  $B_{\mathcal T}(\mathcal{C})$ est la classe des bornes de
$\mathcal{C}$ alors   $\mathcal{C} = Forb_{\mathcal T}(B_{\mathcal T}(\mathcal{C}))$.  Si $\mathcal C$ est une  classe h\'er\'editaire de tournois finis incluse dans $\mathcal F_m$ alors,  d'apr\`es le test obtenu en \cite{pouzet2} et le th\'eor\`eme de Higman sur les mots, ses bornes, compt\'ees \`a des isomorphismes pr\`es, sont en nombre fini. En particulier:
\begin{theoreme} \label{THM fini de bornes}
Les bornes de la classe  $\mathcal I_{m}^{<\omega}$, compt\'ees
\`a des isomorphismes pr\`es, sont en nombre fini.
\end{theoreme}

Le 3-cycle $C_{3}$ est, \`a des isomorphismes pr\`es, l'unique borne de  $\mathcal{I}_{0}^{{< \omega}}$ dans $\mathcal T$. Soient $D_{5} := (\underline{5}, (\{1,3\},\{0,4\}))$ et $C_{3}.\underline {2} := (\underline{6}, (\{0,2\},\{3,5\}))$.
Soit $P_{7}$ le tournoi de Paley d\'efini sur $\mathbb{Z}/7\mathbb{Z}$ par
$A(P_{7}) := \{(i,j): j-i \in \{1,2,4\}\}$ et soit $B_{6} := P_{7} - 6$.
Observons que $B_{6}$ est isomorphe \`a $Inv(\underline{6}, (\{0,3,4\}, \{1,4,5\}))$. Le th\'eor\`eme suivant d\'ecoule du th\'eor\`eme de d\'ecomposition de Gallai \cite{Galai} pour les tournois et d'un r\'esultat de Latka \cite{Latka} caract\'erisant les tournois finis, ind\'ecomposables et n'abritant pas $V_{5}$.

\begin{theoreme} \label{THM 5 bornes}
\`A des isomorphismes pr\`es, les bornes de la classe $\mathcal {I}_{1}^{<\omega}$ sont les tournois $B_{6}$, $C_{3}.\underline{2}$, $D_{5}$ et les tournois critiques $T_{5}$ et $V_{5}$.  En particulier, $\mathcal I_{1}^{<\omega} = Forb_{\mathcal T}(\{B_{6}, C_{3}.\underline {2}, D_{5}, T_{5}, V_{5}\})$.
\end{theoreme}

Un tournoi  est  \emph{acycliquement  ind\'ecomposable} s'il n'a pas d'intervalle acyclique  ayant plus d'un \'el\'ement \cite{C-J}.
Tout tournoi est une somme lexicographique de tournois acycliques  index\'ee par un tournoi acycliquement  ind\'ecomposable \cite{B-P}.
Donc, en vertu du lemme \ref{sum},   \emph{les \'el\'ements  de  $\mathcal{I}_{m}$ sont les sommes lexicographiques de tournois acycliques index\'ees par des tournois acycliquement  ind\'ecomposables appartenant \`a $\mathcal {I}_{m}$}. De fa\c{c}on \'equivalente, {\it les bornes  de  $\mathcal{I}_{m}^{< \omega}$ dans  $\mathcal T$ sont acycliquement ind\'ecomposables}.

\begin{theoreme}
Un tournoi  $T$ ayant au moins deux sommets est un \'el\'ement  acycliquement ind\'ecomposable de $\mathcal {I}_{1}^{<\omega}$ si et seulement si $T$ est isomorphe \`a $U_{2n+1}$, $\underline{2}(\underline{1}, U_{2n+1})$, $\underline{2}(U_{2n+1}, \underline{1})$ ou $\underline{3}(\underline{1}, U_{2n+1}, \underline{1})$, o\`u $n \geq 1$ et $U_{2n+1}: = Inv(\underline{2n+1}, 2\mathbb{N}_{<n+1})$.
\end{theoreme}

La classe $\mathcal O^{<\omega}_{m}$ est une classe de Fra\"{\i}ss\'e (c'est-\`a-dire est h\'er\'editaire et a la propri\'et\'e d'amalgamation) donc $\mathcal O^{\leq \omega}_{m}$ contient une structure d\'enombrable homog\`ene et d'\^age $\mathcal O^{<\omega}_{m}$.  Cette structure \'etant unique \`a l'isomorphie pr\`es, notons  la  $C(m)$. Soit $W(m):=Inv(C(m))$. Le tournoi $W(m)$ est un tournoi d\'enombrable d'indice $m$ qui abrite tous les tournois de la classe $\mathcal{I}_{m}^{\leq \omega}$. Si  $m= 0$, $W(0)$ est $\underline{\mathbb{Q}}$, la cha\^{\i}ne des nombres rationnels.
Si $m \geq 1$, consid\'erons $m$ nombres r\'eels
$\alpha_{0}, \cdots, \alpha_{m-1}$ tels que $1, \alpha_{0}, \cdots, \alpha_{m-1}$ soient {\it
rationnellement ind\'ependants}, c'est-\`a-dire tels que pour tous nombres
rationnels $\beta, \beta_{0}, \cdots, \beta_{m-1}$, si $\beta =
\displaystyle \sum_{i<m} \beta_{i} \alpha_{i}$, alors pour
tout $i <m$, $\beta = \beta_{i} = 0$. Soit  $f \in
2^{\mathbb{N}_{<m}}$ (l'ensemble des applications de $\mathbb{N}_{<m}$ dans $\{0, 1\}$). Posons  $\alpha(f) := \sum_{i<m} f(i) \alpha_{i}$, $\mathbb{Q}_{f} := \mathbb{Q}
+ \alpha(f)$, $\mathbb{Q}(m) :=
\bigcup_{f \in 2^{\mathbb{N}_{<m}}}\mathbb{Q}_{f}$,
o\`u pour $a \in \mathbb{R}$, $\mathbb{Q}
+ a := \{r+a: r \in \mathbb{Q}\}$ et soit $\underline{\mathbb{Q}}(m)$ le tournoi acyclique induit sur $\mathbb{Q}(m)$ par l'ordre naturel sur les r\'eels. Pour tout $i<m$, posons  $X_{i} := \bigcup\limits_{\{f \in 2^{\mathbb{N}_{<m}} : f(i) = 1\}}\mathbb{Q}_{f}$.  Alors $C(m)=(\underline{\mathbb{Q}}(m), (X_{i})_{i<m} )$.

\section*{Remerciements}Les auteurs remercient A. Bondy et S. Thomass\'e pour leur soutien et leurs suggestions. Ils remercient l'arbitre pour son examen tr\`es fouill\'e, la correction des inexactitudes, ses commentaires et suggestions.

\end{document}